\documentclass[runningheads]{llncs}

\usepackage[T2A]{fontenc}
\usepackage[english]{babel}

\usepackage{url}
\usepackage{multirow}
\usepackage{amssymb}
\usepackage{amsxtra}
\usepackage{makeidx}
\usepackage{amsfonts}
\usepackage{mathtools}
\usepackage{algorithm}
\usepackage{algorithmic}
\usepackage{float}
\usepackage{graphicx}
\usepackage{comment}
\usepackage[font=small,labelfont=bf]{caption}
\usepackage{subcaption}
\usepackage{breqn}

\usepackage{color}
\usepackage{hyperref}

\DeclareUnicodeCharacter{202A}{\,}

\begin{document}
	
\title{Gradient-Type Methods for Optimization Problems with Polyak-Łojasiewicz Condition: Early Stopping and Adaptivity to Inexactness Parameter  
\thanks{The research was supported by Russian Science Foundation and Moscow (project No. 22-21-20065, https://rscf.ru/project/22-21-20065).}}
\titlerunning{Gradient-Type Adaptive Methods}

    \author{Ilya A. Kuruzov \inst{1}\orcidID{0000-0002-2715-5489} \and Fedor S. Stonyakin  \inst{1,2} \orcidID{0000-0002-9250-4438} \and Mohammad S. Alkousa \inst{1,3}\orcidID{0000-0001-5470-0182}}
    
\authorrunning{I. Kuruzov et al.}

\institute{Moscow Institute of Physics and Technology, Moscow, Russia
	\and
	V.\,I.\,Vernadsky Crimean Federal University, Simferopol, Russia
	\and 
	HSE University, Moscow, Russia
	\\
	\email{kuruzov.ia@phystech.edu,  fedyor@mail.ru, mohammad.alkousa@phystech.edu}
	}

\maketitle
\begin{abstract}
Due to its applications in many different places in machine learning and other connected engineering applications, the problem of minimization of a smooth function that satisfies the Polyak-Łojasiewicz condition receives much attention from researchers. Recently, for this problem, the authors of \cite{polyakstopping2022} proposed an adaptive gradient-type method using an inexact gradient. The adaptivity took place only with respect to the Lipschitz constant of the gradient. In this paper, for problems with the Polyak-Łojasiewicz condition, we propose a full adaptive algorithm, which means that the adaptivity takes place with respect to the Lipschitz constant of the gradient and the level of the noise in the gradient. We provide a detailed analysis of the convergence of the proposed algorithm and an estimation of the distance from the starting point to the output point of the algorithm. Numerical experiments and comparisons are presented to illustrate the advantages of the proposed algorithm in some examples.


\end{abstract}

\section{Introduction}

With the increase in the number of applications that can be modeled as large- or even huge-scale optimization problems (some of such applications arise in machine learning, deep learning, data science, control, signal processing, statistics, and so on), first-order methods, which require low iteration cost as well as low memory storage, have received much interest over the past few decades \cite{beck2017}. Gradient-type methods may be regarded as the cornerstone and core of the numerical methods for solving optimization problems. 

For the problem of  minimization of a smooth function $f$, it is well known that if $f$ is strongly-convex, then the gradient descent method achieves a global linear convergence rate  \cite{nesterovbook}. However, many of the fundamental models in machine learning such as least squares and logistic regression yield objective functions that are convex but not strongly convex. This matter led to the formation of motivation for seeking some alternatives to strong convexity, showing that it is possible to obtain linear convergence rates for problems such as least squares and logistic regression. One of these alternatives is the Polyak-Łojasiewicz inequality (or PL-condition). This inequality was originally introduced by B.~Polyak \cite{Polyak_1963}, who proved that it is sufficient to show the global linear convergence rate for the gradient descent without assuming convexity.  
PL-condition is very well studied by many researchers in many different works for many different settings of optimization problems and has been theoretically verified for objective functions of optimization problems arising in many practical problems. For example, it has been proven to be true for objectives of over-parameterized deep networks \cite{deepneural2018}, learning LQR models \cite{fazelglobal2018}, and phase retrieval \cite{sungeometric2018}. More discussions of PL-condition and many other simple problems can be found in  \cite{karimilinear2016}. Note that many other important classes of non-convex problems (Lipschitz problems, weakly convex problems, weakly $\alpha$-quasiconvex functions) are investigated by different authors (see e.g. \cite{sergeyev2020safe,strongin2013global,sungeometric2018}).

In the first-order methods (thus the  classical gradient descent method), the availability of an exact first-order oracle is assumed. That is, the oracle must provide at each given point the exact values of the function and its gradient. But unfortunately, in many applications, there is no access to this exact information (especially information about the gradient) at each iteration of the method. This led researchers to investigate the behavior of first-order methods which have the possibility to work with an inexact oracle. This problem attracted the attention of many researchers in mathematical optimization. In \cite{Devolderfirst2014} (which can be considered as a fundamental work in this direction) the authors introduce the notion of an inexact first-order oracle, which naturally appears in the context of smoothing techniques, Augmented Lagrangians, and many other situations. See \cite{Devolderfirst2014,Devolderthesis,Aspremont2008,kuruzov2021sequential,Stonyakininexact2020,Vasin2021} and references therein, for more details. 

It is known that the analysis of the convergence of the gradient descent method, implied the constant step-size which depends on the Lipschitz constant of the gradient of the objective function (constant of smoothness). But in many applied optimization problems, it is difficult to estimate this constant. For example, the well-known Rosenbrock function and its multidimensional generalizations (for example, the Nesterov-Skokov function \cite{Nest_Skok}) have only a locally Lipschitz continuous gradient. Thus, it is impossible to estimate the Lipschitz constant of the gradient for these functions without additional restrictions on the operation region of the method. In order to overcome the difficulty to estimate the value of the Lipschitz constant of the gradient, there have been many attempts to construct a method with adaptivity on the step-size, see for example \cite{gasnikovbook2018,sergeyev2020safe,Stonyakininexact2020,strongin2013global}.

Recently, for the problem of minimizing a smooth function that satisfies the PL-condition (which is the problem under consideration in this paper), in \cite{polyakstopping2022} the authors proposed non-adaptive and adaptive gradient-type methods using the notion of the non-exact gradient. They analyzed the proposed algorithms and the influence of non-exactness on the rate of convergence. But the adaptivity takes place with respect to the Lipschitz constant of the gradient, where it is still necessary to calculate the level of the noise in the gradient exactly. 

In this paper, we continue the research in order to construct an adaptive gradient-type method that was studied in \cite{polyakstopping2022} for the first time and propose an adaptive algorithm for problems with objective functions that satisfy  PL condition with a detailed analysis of its convergence and an estimation of the distance from the starting point to the output point of the algorithm. The adaptivity in the proposed algorithm in this paper will be in both parameters: the Lipschitz constant of the gradient and the level of the noise in the gradient. Therefore, the proposed algorithm is fully adaptive.

This paper consists of an introduction and 4 main sections. In Sect. \ref{basic} we formulate basic concepts, definitions, and assumptions that are connected with the problem under consideration. In Sect. \ref{halfadptive} we mention an adaptive algorithm that was proposed in \cite{polyakstopping2022}, the adaptivity in this algorithm is for the Lipschitz constant of the gradient of the objective function. Sect. \ref{fulladaptive} is devoted to a fully adaptive algorithm (the adaptivity in the Lipschitz constant of the gradient and the level of the noise in the gradient) for problems with objective functions that satisfy PL-condition. The last Sect. \ref{experiments_Alg5} is devoted to the numerical experiments which demonstrate the effectiveness of the proposed algorithm. The conducted experiments were conducted for the minimization problem of the quadratic form, the logistic regression problem, and for one minimization problem connected with the system of nonlinear equations.


\section{Problem Statement and Basic definitions}\label{basic}

In this section, we present the problem statement and some basic concepts and definitions.

\begin{definition}\label{def_Lsmooth}
The differentiable function $f: \mathbb{R}^n \longrightarrow   \mathbb{R}$ is an $L$-smooth (or $\nabla f$ is Lipschitz-continuous) w.r.t. $\|\cdot\|$, for some constant $L>0$, if
\begin{equation}\label{f2}
\|\nabla f(x)-\nabla f(y)\|\leqslant L\|x-y\|\quad \forall\,x,\,y\in\mathbb{R}^{n}.
\end{equation}
Here and everywhere in the paper, the norm $\|\cdot\|$ indicates the Euclidean norm.
\end{definition}

\begin{definition}\label{def_PLcond}
Let $f$ be an $L$-smooth function. The gradient $\nabla f$ satisfies the Polyak-Łojasiewicz condition (for brevity, we write the PL-condition) \cite{Polyak_1963} if the following inequality holds 
\begin{equation}\label{f1}
f(x)-f^{*}\leqslant\frac{1}{2\mu}\|\nabla f(x)\|^{2}\quad \forall\,x\in\mathbb{R}^{n},
\end{equation}
where $\mu > 0, f^* = f(x_*)$ and $x_*$ is one of the exact solutions of the optimization problem under consideration.
\end{definition}

In this paper, we will consider a classical optimization problem 
\begin{equation}\label{main_problem}
    \min_{x \in \mathbb{R}^n} f(x),
\end{equation}
when the objective function $f$ satisfies \eqref{f2} and \eqref{f1}.

Let us denote by $\widetilde{\nabla}f(x)$ the approximate value of $\nabla f(x)$ at any requested point $x$ (alternatively we also call $\widetilde{\nabla}f(x)$ an inexact gradient of $f$ at $x$), this means 
\begin{equation}\label{inexact_grad}
    \nabla f(x)=\widetilde{\nabla}f(x)+ v(x),\quad\text{and}\quad \|v(x)\|\leqslant\Delta,
\end{equation}
for some fixed $\Delta>0$. 
From this, \eqref{f1} means
\begin{equation}\label{f3}
f(x)-f^{*}\leqslant\frac{1}{\mu}(\|\widetilde{\nabla}f(x)\|^{2}+\Delta^{2})\quad \forall\,x\in\mathbb{R}^{n}.
\end{equation}

\section{Gradient Descent with an Adaptive Step-Size Policy}\label{halfadptive}

Due to the difficulty of estimating the Lipschitz constant of the gradient of the objective function, the researchers are actively working in order to propose methods that overcome this difficulty in adaptive forms. In \cite{polyakstopping2022}, the authors proposed an adaptive algorithm (listed as algorithm \ref{adapt_gd}, below), which is a generalization of the universal gradient method  \cite{NesterovUniversal} for working with an inexact gradient of the  functions satisfying the PL-condition. The adaptivity in Algorithm \ref{adapt_gd}, is for the Lipschitz constant $L$. 

For the inexact gradient \eqref{inexact_grad} we can get an inequality similar to \eqref{f2}, as follows:
\begin{align*}
    f(x)\leqslant f(y) + \langle \widetilde{\nabla} f(y), x-y\rangle + L\|x-y\|^2 +\frac{\Delta^2}{2L}, \quad  \forall x,y\in\mathbb{R}^n.
\end{align*}

This inequality contains an exact calculation of the value of the function $f$ at an arbitrary point from the domain of definition. 

Let us assume that we can calculate the inexact value $\widetilde{f}$ of the function $f$ at any point $x$, so that
\begin{equation}
    \label{inexactf_cond}
    |f(x)-\widetilde{f}(x)|\leqslant \delta,
\end{equation}
for some $\delta >0$. Then the following inequality holds:
\begin{equation}
    \label{smooth_cond_Lk}
    \widetilde{f}(x)\leqslant \widetilde{f}(y)+\langle\widetilde{\nabla} f(y), x-y\rangle + L\|x-y\|^2 + \frac{\Delta ^2}{2L} + 2\delta,   \;\;\; \forall x,y\in\mathbb{R}^n.
\end{equation}


\begin{algorithm}[H]
\caption{Adaptive Gradient Descent with Inexact Gradient \cite{polyakstopping2022}.}
\label{adapt_gd}
  \begin{algorithmic}[1]
  \REQUIRE {$x_0$, } $L_{\min} \geqslant 0, L_0 \geqslant L_{\min}, \delta \geqslant 0, \Delta \geqslant 0$.
  \STATE Set $k:=0$
  \STATE Calculate 
  \begin{equation}
  \label{x_k}
  x_{k+1} = x_k - \frac{1}{2L_k}\widetilde{\nabla}f(x_k).
  \end{equation}
    \STATE If 
      \begin{equation}
      \label{cond_it}
    \widetilde{f}(x_{k+1})\leqslant\widetilde{f}(x_k)+\langle\widetilde{\nabla} f(x_k), x_{k+1}-x_k\rangle + L_k\|x_{k+1}-x_k\|^2 + \frac{\Delta^2}{2L_k} +2\delta, \end{equation}
    then $k:=k+1$, $L_k := \max\left\{\frac{L_{k-1}}{2}, L_{\min}\right\}$ and go to Step 2. Otherwise, set $L_k:=2L_k$ and go to Step 3.
\RETURN $x_k$.
\end{algorithmic}
\end{algorithm}

For Algorithm \ref{adapt_gd}, with a sufficiently small inexact gradient, at each point in the sequence $\{x_k\}_{k \geqslant 0}$, 
\begin{equation}
\label{stop_cond_adaptGD}
    \|\widetilde{\nabla}f(x_k)\|\leqslant 2\Delta,
\end{equation}
and according to \eqref{f3}, we can guarantee that $f(x_k)-f^*\leqslant \frac{5\Delta^2}{\mu}$.


\bigskip

For Algorithm \ref{adapt_gd}, in \cite{polyakstopping2022}, the following theorem has been proved.
\begin{theorem}\label{theorem:adaptL_inexactf}
Suppose  that $f$ satisfies PL-condition \eqref{f1} and conditions \eqref{inexactf_cond}, $\Delta^2 \geqslant16L\delta$ hold. Let  the parameter $L_{\min}$ in Algorithm \ref{adapt_gd} be such that $L_{\min}\geqslant \frac{\mu}{4}$ and one of the following holds:
\begin{enumerate}
    \item Algorithm \ref{adapt_gd} works $N_*$ steps where $N_*$ is such that
\begin{equation}
    \label{N_star_adapt_gd}
    N_*= \left\lceil\frac{8L}{\mu}\log \left(\frac{\mu (f(x_0)-f^*)}{\Delta^2}\right)\right\rceil.
\end{equation}
\item For some $N \leqslant N_*$, at the $N$-th iteration of Algorithm \ref{adapt_gd},  stopping criterion \eqref{stop_cond_adaptGD} is satisfied for the first time.
\end{enumerate}
Then for the output point $\widehat{x}$ ($\widehat{x} = x_{N}$ or $\widehat{x} = x_{N_*}$) of Algorithm \ref{adapt_gd} will satisfy the following inequalities
$$
f(\widehat{x})-f^{*}\leqslant\frac{5\Delta^{2}}{\mu},
$$
\begin{equation}
\label{dist}
    \|\widehat{x}-x_0\|\leqslant\frac{8\Delta}{\mu} \sqrt{\frac{\gamma^2}{2}  + \frac{4\gamma L}{\mu }} \log \left(\frac{\mu (f(x_0)-f^*)}{\Delta^2}\right) + \frac{16\sqrt{\gamma L(f(x_0)-f^*)}}{\mu},
\end{equation}
where $\gamma=\frac{L}{L_{\min}}$. Also, the total number of calls to the subroutine for calculating inexact values of the objective function and \eqref{x_k} is not more than $2N+\log\frac{2L}{L_0}$.
\end{theorem}

Note, that Algorithm \ref{adapt_gd} uses subroutines for finding the inexact value of the objective function more often than the gradient method with a constant step-size. But the number of calls to these subroutines in Algorithm \ref{adapt_gd} is not more than  $2N+\log\frac{2L}{L_0}$. This means that the "cost" of an iteration of the adaptive Algorithm \ref{adapt_gd} is, on average, comparable to about two iterations of the non-adaptive method (i.e. with a constant step-size). At the same time, the accuracy achieved by Algorithm \ref{adapt_gd} and the non-adaptive one also equals approximately \cite{polyakstopping2022}.

\section{Gradient Descent with Adaptivity in the Step-Size and Inexactness of the Noise Level}\label{fulladaptive}

As we saw in the previous section, the adaptivity in Algorithm \ref{adapt_gd} was for the parameter $L$ only. In this section,  we consider an Algorithm (listed as Algorithm \ref{adapt_gd_delta}, below), which is a generalization of  Algorithm \ref{adapt_gd} for the case of an unknown noise level $\Delta_k$. It means that the adaptivity in Algorithm \ref{adapt_gd_delta}, will be for two parameters $L$ and $\Delta$. In this algorithm, at each iteration, we select the constants $\Delta_k$ and $L_k$, at each iteration, such that they satisfy the inequality for smooth functions with an inexact gradient for points from neighboring iterations.

\begin{algorithm}[H]
	\caption{Adaptive Gradient Descent Method for unknown $L$ and $\Delta$.}
   \label{adapt_gd_delta}
  \begin{algorithmic}[1]
  \REQUIRE {$x_0$,} $L_{\min} \geqslant \frac{\mu}{4} > 0, L_0 \geqslant L_{\min}, \Delta_0 > 0,  \Delta_{\min} > 0$.
  \STATE {Set} $k:=0$.
  \STATE {Set} $L_k := \max\left\{\frac{L_{k-1}}{2}, L_{\min}\right\}$.
  \STATE Calculate:
  \begin{equation*}
  x_{k+1} = x_k - \frac{1}{2L_k}\tilde{\nabla}f(x_k).
  \end{equation*}
    \STATE If
    \begin{equation}\label{cond_it_L_delta}
        f(x_{k+1})\leqslant f(x_k) + \langle\tilde{\nabla} f(x_k), x_{k+1}-x_k\rangle + \Delta_k \|x_{k+1}-x_k\| + \frac{L_k}{2}\|x_{k+1}-x_k\|^2,
    \end{equation}
    then go to Step 5. Otherwise, set $L_k:=2L_k$ and $\Delta_k = 2\Delta_k$ and go to step 3. 
    \STATE Find the minimal value of $\Delta_k$, such that \eqref{cond_it_L_delta} holds, and $\Delta_k\geq\Delta_{\min}$, also $\Delta_k\geqslant\max_{j<k}\Delta_j$, for $k\geqslant1$.
    \STATE $L_k:=\max\left(\frac{L_k}{2}, L_{\min}\right); x_{k+1} = x_k - \frac{1}{2L_k}\tilde{\nabla}f(x_k)$.
    \STATE If  \eqref{cond_it_L_delta} holds and $L_k \geqslant L_{\min}$, then  go to Step 6.
    Otherwise, set $k:=k+1$ and go to Step 3.
\RETURN $x_k$.
\end{algorithmic}
\end{algorithm}
 
For the sequence of points generated by Algorithm \ref{adapt_gd_delta}, due to fulfillment of condition \eqref{cond_it_L_delta}, the following inequality holds
\begin{equation} \label{k_inequality_delta}
    f(x_{k+1})-f(x_k)\leqslant\frac{\Delta_k^2}{2L_k} - \frac{1}{4L_k}\|\tilde{\nabla}f(x_k)\|^2.
\end{equation}

In this case, using \eqref{k_inequality_delta}, we can get the following estimate
\begin{equation}\label{inexact_fk_delta}
    f(x_{k+1})-f^*\leqslant \prod_{j=0}^k \left(1-\frac{\mu}{4L_j}\right)(f(x_0)-f^*) + \frac{\max_{j \leqslant k}\Delta_{j}^2}{\mu},
\end{equation}
or by $L_j \leqslant \max_{j \leqslant k} L_j$, we get
\begin{equation}\label{inexact_fk_delta_maxL}
    f(x_{k+1})-f^*\leqslant  \left(1-\frac{\mu}{4\max_{j\leqslant k} L_j}\right)^{k+1}(f(x_0)-f^*) + \frac{\max_{j\leqslant k}\Delta_{j}^2}{\mu}.
\end{equation}

Here and below, by $\max_{j\leqslant k}\Delta_{j}^2$ we mean the maximum of all estimates $\Delta_j$ up to the $k$-th iteration.

Let us estimate the value $\max_{j\leqslant k} L_j$. Note that inequality \eqref{cond_it_L_delta} holds for $L_k\geqslant L$ and for $\Delta_k\geqslant \Delta$ for all $k \geqslant 0$.

Let us consider an arbitrary $K$-th iteration. If $\frac{\Delta}{\Delta_K}\leqslant\frac{L}{L_K},$ then when $L_K$ reaches $L$, $\Delta_K$ reaches $\Delta$ also, it will not be more iterations in the step 4. Hence, if $\frac{\Delta}{\Delta_K}\leqslant\frac{L}{L_K},$ then $L_K\leqslant2L$. On the other hand, if $\frac{\Delta}{\Delta_K}> \frac{L}{L_K}$, then in the worst case $L_K = L$ at the beginning and after the completion of the process, we have $L_K\leqslant\frac{2 \Delta}{\Delta_K}L\leqslant\frac{2\Delta}{\Delta_{\min}}L$. Thus, we obtain the estimate $L_{\max} \leqslant 2L\max\left\{\frac{\Delta}{\Delta_{\min}},1\right\}$. Denoting $ L_{\max} = L\max\left\{\frac{\Delta}{\Delta_{\min}}, 1\right\}$, we obtain the following refinement of the bound \eqref{inexact_fk_delta_maxL}:
\begin{equation} \label{inexact_fk_delta_maxL2}
    f(x_{k+1})-f^*\leqslant   \left(1-\frac{\mu}{8L_{\max}}\right)^{k+1}\left(f(x_0)-f^*\right) + \frac{\max_{j\leqslant k}\Delta_{j}^2}{\mu}.
\end{equation}

In a similar way, for $\Delta_k$ we obtain that  $\Delta_k \leqslant \Delta_{\max}:= 2\Delta \max \left\{\frac{L}{L_{\min}}, 1 \right\}$. Thus, $\max_{j\leqslant k} \Delta_j \leqslant \Delta_{\max}$.

In this case, we can stop the method after reaching the accuracy by the gradient $\|\tilde{\nabla}f(x_k)\|\leqslant 2\max_{j \leqslant k} \Delta_j$, which guarantees an estimate for the accuracy by the function 
$$
    f(x_k)-f^*\leqslant  \frac{5\max_{j\leqslant k} \Delta_j^2}{\mu}\leqslant  \frac{5\Delta_{\max}^2}{\mu}.
$$ 

On the other hand, using the estimate \eqref{inexact_fk_delta_maxL2} and the introduced notation, we can guarantee the following rate of convergence
\begin{equation}\label{inexact_fk_delta_maxL3}
    f(x_{k+1})-f^*\leqslant  \left(1-\frac{\mu}{8L_{\max}}\right)^{k+1}\left(f(x_0)-f^*\right) + \frac{\Delta_{\max}^2}{\mu}.
\end{equation}
Then we get the following expression for the number of iterations
$$
    N \leqslant \left\lceil \frac{8L_{\max}}{\mu} \log \left(\frac{\mu (f(x_0)-f^*)}{4\Delta_{\max}^2}\right) \right\rceil.
$$

As before, we get an estimate for the distance from the starting point to the current one, we can estimate it as follows:
$$
    \|x_N-x_0\| \leqslant \frac{8 \Delta_{\max}}{\mu}  \sqrt{\frac{\gamma^2}{2}  + \frac{2\gamma L_{\max}}{\mu }} \log \left(\frac{\mu (f(x_0)-f^*)}{4\Delta_{\max}^2}\right) + \frac{16\sqrt{\gamma L_{\max}(f(x_0)-f^*)}}{\mu},
$$
where $\gamma=\frac{L}{L_{\min}}$.

We can also estimate the total number of evaluations of the  function $f$ at each iteration of Algorithm \ref{adapt_gd_delta}. As mentioned earlier, the condition  \eqref{cond_it_L_delta} is satisfied if $\Delta_k \geqslant \Delta$ and $L_k \geqslant L$. Thus, step 4 will be repeated at one iteration no more than 
$$
    1 + \log_2\left(\max\left\{ \frac{L}{L_{\min}},\frac{\Delta}{\Delta_{\min}}\right\}\right).
$$
Step 5 will be repeated no more than $\log_2 \left(\frac{2L}{L_{\min}}\right)$ times. Thus, the total number of function evaluations does not exceed 
$$
    N_* \log_2 \left(\frac{4L}{L_{\min}}\max\left\{ \frac{L}{L_{\min}},\frac{\Delta}{\Delta_{\min}}\right\}\right).
$$

In this case, we get the following result about the work of Algorithm \ref{adapt_gd_delta}.
\begin{theorem}\label{theorem:adaptLdelta}
Let Algorithm \ref{adapt_gd_delta} works either
$$
    N_*= \left\lceil\frac{8L_{\max}}{\mu}\log \left(\frac{\mu (f(x_0)-f^*)}{4\Delta_{\max}^2}\right)\right\rceil
$$
steps, or for some $k_* \leqslant N_*$ on the $k_*$-th iteration of Algorithm \ref{adapt_gd_delta}, the stopping criterion  $\|\tilde{\nabla}f(x_{k_*})\|\leqslant 2\max_{j\leqslant k} \Delta_j$ be achieved.
Then for the output point $\widehat{x}$ ($\widehat{x} = x_{k_*}$ or $\widehat{x} = x_{N_*}$) of the Algorithm \ref{adapt_gd_delta}, the following inequality is guaranteed to be true
$$
    f(\widehat{x})-f^{*}\leqslant \frac{5\max_{j \leqslant k} \Delta_j^2}{\mu} \leqslant \frac{5\Delta_{\max}^2}{\mu}.
$$

Moreover,
$$
    \|\widehat{x}-x_0\|\leqslant\frac{8\Delta_{\max}}{\mu} \sqrt{\frac{\gamma^2}{2}  + \frac{2\gamma L_{\max}}{\mu }} \log \left(\frac{\mu (f(x_0)-f^*)}{4\Delta_{\max}^2}\right) + \frac{16\sqrt{\gamma L_{\max}(f(x_0)-f^*)}}{\mu},
$$
where  $\gamma=\frac{4L_{\max}}{\mu}$, $\Delta_{\max} = 2 \Delta \max \left\{\frac{L}{L_{\min}}, 1\right\}$, $ L_{\max} = L \max \left\{ \frac{\Delta}{\Delta_{\min}},1\right\}$. Also, the total number of calls to the calculation of the function $f$ is not more than 
$$
    N_* \log_2 \left( \frac{4L}{L_{\min}} \max \left\{ \frac{L}{L_{\min}}, \frac{\Delta}{\Delta_{\min}} \right\} \right).
$$
\end{theorem}

\begin{remark}
The value $L_{\max}$ estimates the maximum value of the parameter $L_k$. The estimates above remain valid if $L_{\max}$ is replaced by $\max_{j\leqslant k}L_j$.
\end{remark}



\begin{remark}
Let at any point $x$ we have a model $(\tilde{f}, \tilde{\nabla}f)$ of the function $f$ such that the following conditions are satisfied:
\begin{equation}\label{smooth_cond}
    \tilde{f}(x)\leqslant \tilde{f}(y)+\langle \tilde{\nabla} f(y), x-y\rangle + \frac{L}{2}\|x-y\|^2 + \Delta \|x-y\| + \delta,
\end{equation}
and
\begin{equation}
    \label{PL_Delta_delta}
    \tilde{f}(x)-f^*\leqslant \frac{1}{\mu}\left(\|\tilde{\nabla} f(x)\| +\Delta^2\right)+\delta.
\end{equation}
Note that if $\tilde{\nabla}f$ is a $\Delta$-inexact gradient and $\tilde{f}$ is a function such that $|f(x)-\tilde{f}(x)|\leqslant \delta$
for each $x$, then the conditions above are satisfied. A natural modification of the Algorithm  \ref{adapt_gd} is the selection of $L$ such that \eqref{smooth_cond} will be satisfied for it. In this case, the achieved accuracy will decrease from $\frac{5\Delta^2}{\mu}$ to $\frac{5\Delta^2}{\mu} +\delta$. We can act similarly in the case of a known constant $\delta$ in Algorithm \ref{adapt_gd_delta}. In this case, the additional factor $\delta$ will also appear exactly. However, if the given parameter $\delta$ is not known, then it can also be selected. However, this will lead to an additional complication of the algorithm, which we will not discuss here. Note that the condition  \eqref{cond_it_L_delta} will be hold if $L_k \geqslant L$ and $\Delta_k \geqslant \Delta +  \frac{2L_k}{\|\tilde{\nabla}f(x_k)\|} \delta $. 
In order for this condition to be achieved, it is enough to modify the update of $\Delta_k$ at the 4th step of Algorithm \ref{adapt_gd_delta} so that $\Delta_k:=2^{\tau}\Delta_k$ for $\tau>1.$ In this case, the estimates for the number of iterations and the accuracy of the solution with respect to the function, and the distance from the starting point to the output point from  Theorem \ref{theorem:adaptLdelta} will remain true, provided that the parameters $\Delta_{\max}$ and $L_{\max} $ will change accordingly.\footnote{In the worst case, $L_k$ at the beginning of a new iteration is already equal to $L$. Then denote by $I_1$ the minimal solution in $I$ of the inequality $2^{\tau I} \Delta_{\min} \geqslant \Delta +  \frac{\delta \sqrt{2}L}{\Delta} 2^I $ provided that $I \geqslant 1$. Then $L_{\max}=2^{I} L$. Similarly, we can get that
$$
    \Delta_{\max} = 2 \left(\Delta + \frac{\delta \sqrt{2} L_{\max}}{\Delta}\right)\cdot \max \left\{ \left(\frac{L}{L_{\min}}\right)^{\tau}, 1\right\}. 
$$}
\end{remark}

\begin{remark}
Note, that the achieved $\Delta_{max}$ can be significantly more than $\Delta$ according to results of Theorem \ref{theorem:adaptLdelta}. Nevertheless, note that we do not research in this work the influence of step 5 of Algorithm \ref{adapt_gd_delta}. Moreover, according to our experiments, the method stops when $\Delta_{\max} \sim \Delta$.
\end{remark}

As before, these estimates give theoretical guarantees for the convergence of our method. But we observe that methods work significantly better in practice than in theory. Particularly, we see in all our experiments that proposed adaptive method \ref{adapt_gd_delta} can approach quality $O\left(\Delta^2\right)$ in terms of gradient norm like gradient method with constant step (see \cite{polyakstopping2022}).

\section{Numerical Experiments}\label{experiments_Alg5}

In this section, in order to demonstrate the performance of Algorithms \ref{adapt_gd}, \ref{adapt_gd_delta} and the algorithm with the constant step-size (see (2.1), (2.9) and Theorems 2, 3 in  \cite{polyakstopping2022}) and compare them, we consider some numerical experiments concerning the quadratic optimization problem,  logistic regression minimization and the solution of the system of nonlinear equations. 

In all experiments, we will use an exact gradient with random noise $v(x)$ in \eqref{inexact_grad}, that is randomly generated, on the $n$ dimensional sphere with radius 1 at the center 0, i.e. $v(x)\sim\mathcal{U}\left(S_1^n(0)\right)$.

All experiments were implemented in Python 3.4, on a computer fitted with Intel(R) Core(TM) i5-8250U CPU @ 1.60GHz, 8000 Mhz, 4 Core(s), 8 Logical Processor(s). The RAM of the computer is 8 GB.

Note, in the work \cite{polyakstopping2022} authors propose the stopping rule $\|\tilde{\nabla}f(x_k)\|\leq \sqrt{6}\Delta$ for known inexactness value $\Delta$. In this work constant $\sqrt{6}$ was chosen for a single criterion of stopping methods and we will use it too.

\subsection{The minimization problem of the quadratic form}

In this subsection, we consider the minimization problem of the quadratic form
\begin{equation}\label{quad_form}
    f(x)=\sum\limits_{i=1}^n d_i x_i^2, \quad x=(x_1,\ldots, x_n) \in \mathbb{R}^n, \, d_i \in \mathbb{R}.
\end{equation}

We run Algorithms \ref{adapt_gd}, \ref{adapt_gd_delta} and the variant with a constant step-size  (we denote this variant ''Alg. constant'' in the listed tables below), for $n=100, \, L=\max_{1 \leqslant i \leqslant n} d_i = 1$ and different values of the parameter $\mu$. 

We take $x_0 =(100,\dots, 100)^\top$ as the initial point of all the compared algorithms.

The results are presented in Tables \ref{tab:qp_NT} and  \ref{tab:qp_dist}.
The results, in Table \ref{tab:qp_NT}  demonstrate the running time (in milliseconds) of the algorithms and the required number of iterations to achieve the accuracy $\|\tilde{\nabla} f(x)\|\leqslant \sqrt{6}\Delta$, for different values of $\Delta$. Meanwhile, the results in Table \ref{tab:qp_dist} demonstrate the achieved accuracy with respect to  $\|\nabla f(x_N)\|$, which is the norm of the gradient of the objective function $f$ at the output point $x_N$ of the algorithms after $N$ iterations, and the distance between the initial point $x_0$ and the output point $x_N$. Note that the distance between $x_0$ and the nearest optimal is equal to $948.7$. 

\begin{table}[htp]
    \centering
    \caption{The results of the algorithms for the quadratic form \eqref{quad_form}, with different values of $\mu$ and $\Delta$, to achieve the accuracy $\|\tilde{\nabla} f(x)\|\leqslant \sqrt{6}\Delta$.}
    \label{tab:qp_NT}
    \begin{tabular}{|c|c||c|c||c|c||c|c|}
      \hline
            &  &\multicolumn{2}{|c||}{Alg. constant \cite{polyakstopping2022}} & \multicolumn{2}{|c||}{Algorithm \ref{adapt_gd}} & \multicolumn{2}{|c|}{Algorithm \ref{adapt_gd_delta}}\\
            \hline
      $\mu$ & $\Delta$ & Iters & Time, ms & Iters &  Time, ms & Iters &  Time, ms \\
      \hline
0.01 & \begin{tabular}{@{}c@{}} $10^{-7}$ \\ $10^{-4}$ \\ $10^{-1}$ \end{tabular}&\begin{tabular}{@{}c@{}} $1525$ \\ $840$ \\ $159$ \end{tabular}&\begin{tabular}{@{}c@{}} $142.31$ \\ $75.99$ \\ $14.81$ \end{tabular}&\begin{tabular}{@{}c@{}} $511$ \\ $301$ \\ $85$ \end{tabular}&\begin{tabular}{@{}c@{}} $226.03$ \\ $158.16$ \\ $33.94$ \end{tabular}&\begin{tabular}{@{}c@{}} $515$ \\ $314$ \\ $170$ \end{tabular}&\begin{tabular}{@{}c@{}} $412.39$ \\ $247.06$ \\ $103.11$ \end{tabular}\\
\hline
0.1 & \begin{tabular}{@{}c@{}} $10^{-7}$ \\ $10^{-4}$ \\ $10^{-1}$ \end{tabular}&\begin{tabular}{@{}c@{}} $169$ \\ $104$ \\ $41$ \end{tabular}&\begin{tabular}{@{}c@{}} $15.96$ \\ $10.15$ \\ $4.46$ \end{tabular}&\begin{tabular}{@{}c@{}} $76$ \\ $49$ \\ $24$ \end{tabular}&\begin{tabular}{@{}c@{}} $29.58$ \\ $18.11$ \\ $8.21$ \end{tabular}&\begin{tabular}{@{}c@{}} $102$ \\ $94$ \\ $54$ \end{tabular}&\begin{tabular}{@{}c@{}} $60.81$ \\ $44.58$ \\ $32.44$ \end{tabular}\\
\hline
0.9 & \begin{tabular}{@{}c@{}} $10^{-7}$ \\ $10^{-4}$ \\ $10^{-1}$ \end{tabular}&\begin{tabular}{@{}c@{}} $11$ \\ $8$ \\ $5$ \end{tabular}&\begin{tabular}{@{}c@{}} $1.57$ \\ $1.28$ \\ $0.90$ \end{tabular}&\begin{tabular}{@{}c@{}} $37$ \\ $26$ \\ $15$ \end{tabular}&\begin{tabular}{@{}c@{}} $15.59$ \\ $19.90$ \\ $7.92$ \end{tabular}&\begin{tabular}{@{}c@{}} $72$ \\ $48$ \\ $39$ \end{tabular}&\begin{tabular}{@{}c@{}} $36.20$ \\ $51.72$ \\ $47.80$ \end{tabular}\\
\hline
0.99 & \begin{tabular}{@{}c@{}} $10^{-7}$ \\ $10^{-4}$ \\ $10^{-1}$ \end{tabular}&\begin{tabular}{@{}c@{}} $6$ \\ $5$ \\ $3$ \end{tabular}&\begin{tabular}{@{}c@{}} $1.02$ \\ $0.89$ \\ $0.52$ \end{tabular}&\begin{tabular}{@{}c@{}} $34$ \\ $24$ \\ $14$ \end{tabular}&\begin{tabular}{@{}c@{}} $13.38$ \\ $9.21$ \\ $5.06$ \end{tabular}&\begin{tabular}{@{}c@{}} $58$ \\ $46$ \\ $48$ \end{tabular}&\begin{tabular}{@{}c@{}} $34.01$ \\ $27.20$ \\ $24.74$ \end{tabular}\\
\hline
\end{tabular}
\end{table}

\begin{table}[htp]
    \centering
    \caption{The results of the algorithms for the quadratic form \eqref{quad_form}, with different values of $\mu$ and $\Delta$. }
    \label{tab:qp_dist}
    \begin{tabular}{|c|c||c|c||c|c||c|c|}
      \hline
            &   &\multicolumn{2}{|c||}{Alg. constant \cite{polyakstopping2022}} & \multicolumn{2}{|c||}{Algorithm \ref{adapt_gd}} & \multicolumn{2}{|c|}{Algorithm \ref{adapt_gd_delta}}\\
            \hline
      $\mu$ &  $\Delta$ & $\|x_N-x_0\|$ & $\frac{\|\nabla f(x_N)\|}{\Delta}$& $\|x_N-x_0\|$ & $\frac{\|\nabla f(x_N)\|}{\Delta}$& $\|x_N-x_0\|$ & $\frac{\|\nabla f(x_N)\|}{\Delta}$\\
      \hline
0.01 & \begin{tabular}{@{}c@{}} $10^{-7}$ \\ $10^{-4}$ \\ $10^{-1}$ \end{tabular}&\begin{tabular}{@{}c@{}} $948.7$ \\ $948.7$ \\ $946.3$ \end{tabular}&\begin{tabular}{@{}c@{}} $2.29$ \\ $2.26$ \\ $2.27$ \end{tabular}&\begin{tabular}{@{}c@{}} $948.7$ \\ $948.7$ \\ $946.7$ \end{tabular}&\begin{tabular}{@{}c@{}} $2.03$ \\ $2.31$ \\ $1.95$ \end{tabular}&\begin{tabular}{@{}c@{}} $948.7$ \\ $948.7$ \\ $947.7$ \end{tabular}&\begin{tabular}{@{}c@{}} $1.84$ \\ $2.43$ \\ $1.07$ \end{tabular}\\
\hline
0.1 & \begin{tabular}{@{}c@{}} $10^{-7}$ \\ $10^{-4}$ \\ $10^{-1}$ \end{tabular}&\begin{tabular}{@{}c@{}} $948.7$ \\ $948.7$ \\ $948.3$ \end{tabular}&\begin{tabular}{@{}c@{}} $2.17$ \\ $2.18$ \\ $2.14$ \end{tabular}&\begin{tabular}{@{}c@{}} $948.7$ \\ $948.7$ \\ $948.3$ \end{tabular}&\begin{tabular}{@{}c@{}} $1.97$ \\ $1.68$ \\ $2.26$ \end{tabular}&\begin{tabular}{@{}c@{}} $948.7$ \\ $948.7$ \\ $948.5$ \end{tabular}&\begin{tabular}{@{}c@{}} $0.87$ \\ $0.80$ \\ $0.83$ \end{tabular}\\
\hline
0.9 & \begin{tabular}{@{}c@{}} $10^{-7}$ \\ $10^{-4}$ \\ $10^{-1}$ \end{tabular}&\begin{tabular}{@{}c@{}} $948.7$ \\ $948.7$ \\ $948.7$ \end{tabular}&\begin{tabular}{@{}c@{}} $0.92$ \\ $0.91$ \\ $0.96$ \end{tabular}&\begin{tabular}{@{}c@{}} $948.7$ \\ $948.7$ \\ $948.6$ \end{tabular}&\begin{tabular}{@{}c@{}} $1.58$ \\ $0.96$ \\ $0.95$ \end{tabular}&\begin{tabular}{@{}c@{}} $948.7$ \\ $948.7$ \\ $948.6$ \end{tabular}&\begin{tabular}{@{}c@{}} $0.69$ \\ $0.79$ \\ $0.72$ \end{tabular}\\
\hline
0.99 & \begin{tabular}{@{}c@{}} $10^{-7}$ \\ $10^{-4}$ \\ $10^{-1}$ \end{tabular}&\begin{tabular}{@{}c@{}} $948.7$ \\ $948.7$ \\ $948.7$ \end{tabular}&\begin{tabular}{@{}c@{}} $0.96$ \\ $0.95$ \\ $1.05$ \end{tabular}&\begin{tabular}{@{}c@{}} $948.7$ \\ $948.7$ \\ $948.7$ \end{tabular}&\begin{tabular}{@{}c@{}} $0.95$ \\ $0.92$ \\ $0.93$ \end{tabular}&\begin{tabular}{@{}c@{}} $948.7$ \\ $948.7$ \\ $948.6$ \end{tabular}&\begin{tabular}{@{}c@{}} $0.71$ \\ $0.68$ \\ $0.69$ \end{tabular}\\
\hline
\end{tabular}
\end{table}

We can see that the adaptive Algorithms \ref{adapt_gd} and \ref{adapt_gd_delta} are slower than the non-adaptive one (algorithm with a constant step-size) for all parameters $\mu$ and $\Delta$. However, they give a gain in the number of iterations for small $\mu$. At the same time, we can notice that the running time for Algorithm \ref{adapt_gd_delta} is longer than for Algorithm  \ref{adapt_gd}. The greatest difference is observed for large values of $\mu$. At the same time, when $\mu$ decreases, they begin to show similar results. Also, we note from Table \ref{tab:qp_dist}, that all three algorithms achieve approximately the same quality, while not going far from the starting point.

\subsection{Logistic Regression}
Now let us examine the work of the compared algorithms in the problem of logistic regression minimization, which has the following form
\begin{equation}\label{logistic_prob}
    \min_{x \in \mathbb{R}^n} f(x) = \frac{1}{m}\sum\limits_{i=1}^m\log \left(1 + \exp\left(-y_i \langle w_i, x\rangle\right)\right ),
\end{equation}
where $y=\left(y_1,\dots, y_m\right)^\top \in [-1,1]^m$ is the feasible variable vector, and  $W=[w_1\dots w_m]\in \mathbb{R}^ {n\times m}$ is the feature matrix, where the vector $w_i\in\mathbb{R}^n$ is from the same space as the optimized weight vector $w$. 

Note that this problem may not have a finite solution in the general case. So we will create an artificial data set such that there is a finite vector $x_*$ minimizing the given function in the way described in \cite{polyakstopping2022}.

\begin{table}[htp]
    \centering
    \caption{The results of Algorithms for the problem \eqref{logistic_prob} with different values of $\Delta$, to achieve the accuracy $\|\tilde{\nabla} f(x)\|\leqslant \sqrt{6}\Delta$.
    }
    \begin{tabular}{|c||c|c|c||c|c|c||c|c|c|}
      \hline
        &\multicolumn{3}{|c||}{Alg. constant \cite{polyakstopping2022}} & \multicolumn{3}{|c||}{Algorithm \ref{adapt_gd}} & \multicolumn{3}{|c|}{Algorithm \ref{adapt_gd_delta}}\\
            \hline
      $\Delta$ & Iters & Time, ms & $\frac{\|\nabla f(\hat{x})\|}{\Delta}$ & Iters &  Time, ms& $\frac{\|\nabla f(\hat{x})\|}{\Delta}$ & Iters &  Time, ms & $\frac{\|\nabla f(\hat{x})\|}{\Delta}$\\
      \hline
\begin{tabular}{@{}c@{}} $10^{-5}$ \\ $10^{-4}$ \\ $10^{-2}$ \end{tabular}&\begin{tabular}{@{}c@{}} $20002$ \\ $9700$ \\ $83$ \end{tabular}&\begin{tabular}{@{}c@{}} $5604.34$ \\ $2605.02$ \\ $49.86$ \end{tabular}&\begin{tabular}{@{}c@{}} $3.56$ \\ $2.42$ \\ $2.29$ \end{tabular}&\begin{tabular}{@{}c@{}} $902$ \\ $472$ \\ $17$ \end{tabular}&\begin{tabular}{@{}c@{}} $2856.98$ \\ $1678.02$ \\ $68.36$ \end{tabular}&\begin{tabular}{@{}c@{}} $2.22$ \\ $2.16$ \\ $2.10$ \end{tabular}&\begin{tabular}{@{}c@{}} $23$ \\ $25$ \\ $17$ \end{tabular}&\begin{tabular}{@{}c@{}} $449.68$ \\ $370.09$ \\ $161.27$ \end{tabular}&\begin{tabular}{@{}c@{}} $3.37$ \\ $3.62$ \\ $0.84$ \end{tabular}\\
\hline
\end{tabular}
\label{tab:log_reg}
\end{table}

We chose $n=200, m=700$, and consider the case of constant inexactness. The results are presented in Table \ref{tab:log_reg}. From these results, we can see that the proposed Algorithm \ref{adapt_gd_delta} stopped faster than gradient descent with constant and adaptive step-size for known $\Delta$. But the obtained quality is a little worse. Nevertheless, we can see the main advantage of the proposed method, that it finds inexactness value without additional information.

\subsection{Solving a System of Nonlinear Equations}
In this subsection, we consider the problem of solving a system of $m$ nonlinear equations
\begin{equation}\label{system_nonlineqs}
    g_i(x)=\sum\limits_{j=1}^n A_{ij} \sin(x_j) + B_{ij} \cos(x_j) := E_i, \quad   i=1, \ldots, m,
\end{equation}
where $x=(x_1,\ldots, x_n) \in \mathbb{R}^n$ and $A_{ij}, B_{ij} \in \mathbb{R}, \; \forall i = 1, \ldots, m; \, j = 1, \ldots, n$. 
This problem can be written as the following optimization problem
\begin{equation}\label{nonlin_pf}
\min_{x\in\mathbb{R}^n} f(x) :=\sum\limits_{i=1}^{m} \left(g_i(x)-E_i\right)^2.
\end{equation}

In the conducted experiments, the parameters $A_{ij}, B_{ij}$ were chosen such that $AB^\top = BA^\top = 0$, where $A, B$ are the matrices with entries $A_{ij}, B_{ij}$, respectively. In this case, the jacobian of the system \eqref{system_nonlineqs} will be $J=A\text{diag}\left(\sin(x)\right)+B\text{diag}\left(\cos(x)\right)$. We can estimate the parameter $\mu$ so that $$
    \mu=\lambda_{\min}(JJ^\top)\geqslant\min\left(\lambda_{\min}(AA^\top),\lambda_{\min}(BB^\top)\right),
$$
where $\lambda_{\min}(X)$  is the smallest eigenvalue of the matrix $X$.

Also, the Lipschitz constant of the gradient of the objective function in \eqref{nonlin_pf}, can be estimated so that $8\sqrt{2}\sigma^2_{\max}\left(\left(A|B\right)\right)$, where $\sigma_{\max}(X)$ is the largest singular value of the matrix $X$. Indeed, let $E \in \mathbb{R}^{m\times n}$ be a matrix formed by $E_i, i = 1, \ldots, m$, then the objective function can be represented as a composition  $f(x)=h(\sin(x), \cos(x))$, where  $h(x,y)=\|Ax+By-E\|^2$, with $L_h=2\sigma^2_{\max}\left(\left(A|B\right)\right)$ as a Lipschitz constant of the gradient $\nabla h$. Then for any $x, y \in \mathbb{R}^n$, we have
\begin{equation*}
    \begin{aligned}
    \|\nabla f(x)-\nabla f(y)\|& \leqslant
    4\|\nabla h(\sin(x), \cos(x))-\nabla h(\sin(y), \cos(y))\| \\ & 
    \leqslant 4L_h \|(\sin(x)-\sin(y), \cos(x)-\cos(y))\| \\&
    \leqslant  4\sqrt{2}L_h\|x - y \|.
\end{aligned}
\end{equation*}

We run the compared algorithms for  $n=256$ and different numbers of equations $m \in \{8,32,128\}$.   We take $x_0 = \mathbf{1}_n=(1,\dots, 1)^\top$ as the initial point of all the compared algorithms. The results are presented in Tables \ref{tab:pf_NT} and  \ref{tab:pf_dist}.
The results, in Table \ref{tab:pf_NT}  demonstrate the running time (in milliseconds) of algorithms and the required number of iterations to achieve the accuracy $\|\tilde{\nabla} f(x)\|\leqslant\sqrt{6}\Delta$, for different values of $\Delta$. Meanwhile, the results in Table \ref{tab:pf_dist}, demonstrate the achieved accuracy with respect to the $\|\nabla f(x_N)\|$, which is the norm of the gradient of the objective function $f$ at the output point $x_N$ of the algorithms after $N$ iterations, and the distance between the initial point $x_0$ and the output point $x_N$.

\begin{table}[htp]
    \centering
    \caption{The results of the algorithms for the problem \eqref{nonlin_pf}, with different values of $\Delta$, to achieve the accuracy $\|\tilde{\nabla} f(x)\|\leqslant \sqrt{6}\Delta$.}
    \begin{tabular}{|c|c|c||c|c||c|c||c|c|}
      \hline
           & &  &\multicolumn{2}{|c||}{Alg. constant \cite{polyakstopping2022}} & \multicolumn{2}{|c||}{Algorithm \ref{adapt_gd}} & \multicolumn{2}{|c|}{Algorithm \ref{adapt_gd_delta}}\\
            \hline
      $m$& $\frac{L}{\mu}$ & $\Delta$ & Iters & Time, ms & Iters &  Time, ms & Iters &  Time, ms\\
      \hline
8 & $2.1\cdot 10^{5}$ & \begin{tabular}{@{}c@{}} $10^{-4}$ \\ $10^{-1}$ \end{tabular}&\begin{tabular}{@{}c@{}} $14434$ \\ $1479$ \end{tabular}&\begin{tabular}{@{}c@{}} $2203.85$ \\ $220.82$ \end{tabular}&\begin{tabular}{@{}c@{}} $508$ \\ $56$ \end{tabular}&\begin{tabular}{@{}c@{}} $447.06$ \\ $45.54$ \end{tabular}&\begin{tabular}{@{}c@{}} $306$ \\ $69$ \end{tabular}&\begin{tabular}{@{}c@{}} $813.99$ \\ $178.63$ \end{tabular}\\
\hline
32 & $5.0\cdot 10^{6}$ & \begin{tabular}{@{}c@{}} $10^{-4}$ \\ $10^{-1}$ \end{tabular}&\begin{tabular}{@{}c@{}} $59643$ \\ $12536$ \end{tabular}&\begin{tabular}{@{}c@{}} $10314.68$ \\ $2290.65$ \end{tabular}&\begin{tabular}{@{}c@{}} $1871$ \\ $383$ \end{tabular}&\begin{tabular}{@{}c@{}} $1808.70$ \\ $393.74$ \end{tabular}&\begin{tabular}{@{}c@{}} $510$ \\ $144$ \end{tabular}&\begin{tabular}{@{}c@{}} $1509.95$ \\ $499.39$ \end{tabular}\\
\hline
128 & $7.7\cdot 10^{8}$ & \begin{tabular}{@{}c@{}} $10^{-4}$ \\ $10^{-1}$ \end{tabular}&\begin{tabular}{@{}c@{}} $921805$ \\ $264361$ \end{tabular}&\begin{tabular}{@{}c@{}} $177759.81$ \\ $47414.07$ \end{tabular}&\begin{tabular}{@{}c@{}} $27405$ \\ $8800$ \end{tabular}&\begin{tabular}{@{}c@{}} $23520.96$ \\ $7367.66$ \end{tabular}&\begin{tabular}{@{}c@{}} $4688$ \\ $1916$ \end{tabular}&\begin{tabular}{@{}c@{}} $11721.18$ \\ $4919.18$ \end{tabular}\\
\hline
\end{tabular}
\label{tab:pf_NT}
\end{table}

\begin{table}[htp]
    \centering
    \caption{The results of the algorithms for the problem \eqref{nonlin_pf}, with different values of $\Delta$.}
    \begin{tabular}{|c|c|c||c|c||c|c||c|c|}
      \hline
            & & &\multicolumn{2}{|c||}{Alg. constant \cite{polyakstopping2022}} & \multicolumn{2}{|c||}{Algorithm \ref{adapt_gd}} & \multicolumn{2}{|c|}{Algorithm \ref{adapt_gd_delta}}\\
            \hline
      $m$ & $\frac{L}{\mu}$ &  $\Delta$ & $\|x_N-x_0\|$ & $\frac{\|\nabla f(x_N)\|}{\Delta}$& $\|x_N-x_0\|$ & $\frac{\|\nabla f(x_N)\|}{\Delta}$& $\|x_N-x_0\|$ & $\frac{\|\nabla f(x_N)\|}{\Delta}$\\
\hline
8 & $2.1\cdot 10^{5}$ & \begin{tabular}{@{}c@{}} $10^{-4}$ \\ $10^{-1}$ \end{tabular}&\begin{tabular}{@{}c@{}} $5.1$ \\ $4.8$ \end{tabular}&\begin{tabular}{@{}c@{}} $2.41$ \\ $2.34$ \end{tabular}&\begin{tabular}{@{}c@{}} $5.1$ \\ $4.8$ \end{tabular}&\begin{tabular}{@{}c@{}} $2.25$ \\ $2.09$ \end{tabular}&\begin{tabular}{@{}c@{}} $5.1$ \\ $4.9$ \end{tabular}&\begin{tabular}{@{}c@{}} $3.07$ \\ $2.28$ \end{tabular}\\
\hline
32 & $5.0\cdot 10^{6}$ & \begin{tabular}{@{}c@{}} $10^{-4}$ \\ $10^{-1}$ \end{tabular}&\begin{tabular}{@{}c@{}} $7.4$ \\ $7.4$ \end{tabular}&\begin{tabular}{@{}c@{}} $2.41$ \\ $2.37$ \end{tabular}&\begin{tabular}{@{}c@{}} $7.5$ \\ $7.4$ \end{tabular}&\begin{tabular}{@{}c@{}} $2.27$ \\ $2.24$ \end{tabular}&\begin{tabular}{@{}c@{}} $7.5$ \\ $7.4$ \end{tabular}&\begin{tabular}{@{}c@{}} $3.31$ \\ $2.72$ \end{tabular}\\
\hline
128 & $7.7\cdot 10^{8}$ & \begin{tabular}{@{}c@{}} $10^{-4}$ \\ $10^{-1}$ \end{tabular}&\begin{tabular}{@{}c@{}} $14.4$ \\ $14.3$ \end{tabular}&\begin{tabular}{@{}c@{}} $2.43$ \\ $2.43$ \end{tabular}&\begin{tabular}{@{}c@{}} $14.4$ \\ $14.3$ \end{tabular}&\begin{tabular}{@{}c@{}} $2.33$ \\ $2.32$ \end{tabular}&\begin{tabular}{@{}c@{}} $14.4$ \\ $14.4$ \end{tabular}&\begin{tabular}{@{}c@{}} $15.67$ \\ $1.84$ \end{tabular}\\
\hline
\end{tabular}
\label{tab:pf_dist}
\end{table}

Now, it can be seen from Table \ref{tab:pf_NT}, that the adaptive Algorithms \ref{adapt_gd} and \ref{adapt_gd_delta} converge much faster than the algorithm with a constant step-size. While the adaptive algorithms converge at approximately the same time, the algorithm with a constant step-size converges at least one and a half times slower. Moreover, for a large number $\frac{L}{\mu}$ it converges more than 10 times slower. A significant efficiency of the adaptive algorithms is observed for a large value of $\frac{L}{\mu}$ (see also the results in Table \ref{tab:pf_dist}).

Note that, we stop the compared algorithms at accuracy $\|\tilde{\nabla} f(x)\|\leqslant \sqrt{6}\Delta$, although estimates for accuracy $\|\tilde{\nabla} f(x)\|\leqslant 2\Delta$ are proved for adaptive algorithms \ref{adapt_gd} (\cite{polyakstopping2022}) and \ref{adapt_gd_delta} (see Theorem \ref{theorem:adaptLdelta}). It was decided to choose a single stopping criterion that both adaptive algorithms can achieve the minimum accuracy from the available ones was chosen. If we consider the criterion $\|\tilde{\nabla}f(x_k)\|\leqslant \sqrt{6}\Delta$, then the results of Theorem \ref{theorem:adaptL_inexactf} and \ref{theorem:adaptLdelta} on the number of iterations and distances from $x_0$ to the output point  $\widehat{x} : = x_N$ remain valid, but they can be refined by increasing the denominator in the logarithm.

\section*{Conclusion}
In this paper, we have considered an adaptive gradient-type method for problems of minimizing a smooth function that satisfies the PL-condition. The adaptivity of the proposed algorithm is in the Lipschitz constant of the gradient and the level of the noise in the gradient. This gives the algorithm the attribute of being fully adaptive. A detailed analysis of its convergence, and an estimation of the distance from the starting point to the output point of the algorithm, were provided. Also, some numerical experiments were conducted for the problem of minimizing a quadratic form, logistic regression, and the problem of solving a system of nonlinear equations.

\end{document}